\documentclass[
a4paper,
eqnright, reqno]{amsart}
\usepackage{graphicx}
\usepackage{lgrind}
\usepackage{hyperref}
\usepackage{myamsart}
\input{mydef}

\newtheorem{Question}[thm]{Question}

\title{Gromov Conjecture on Surface Subgroups:
  Computational Experiments}

\author{Anastasia V. Kisil}
\date{2nd October 2008 }
\address{Trinity College, Cambridge, CB2 1TQ}
\email{ak528@cam.ac.uk}
\thanks{This project was sponsered  by Trinity College, Cambridge.}

\begin{document}
\begin{abstract}
In this paper we investigate Gromov's  question: whether every one-ended
word hyperbolic group contains a surface subgroup.   The case of double
groups is considered by studying the associated
one relator groups. We show that the majority (\(96 \%\)) of the
randomly selected double groups with three generators  have the property. The
experiments are performed on MAGMA software.

\end{abstract}
\AMSMSC{20F67}{20Fxx}
\keywords{one relator groups, surface group, Gromov, word hyperbolic
  group, MAGMA}
\maketitle

\section{Introduction}

In this paper we are going to investigate the following question:
\begin{Question}[Gromov] \cite{problems}
  \label{co:qu}
  Does every one-ended word hyperbolic group contains a surface group?
\end{Question}
Here a ``surface subgroup'' means a subgroup isomorphic to the fundamental
group of a closed surface with non-positive Euler
characteristic. This question of Gromov is
of interest partly because it is a natural generalisation of famous
Surface Subgroup Conjecture. In the case of the fundamental groups of
hyperbolic 3-manifolds it is exactly the Conjecture. 
The question of finding subgroups is studied from
different angles and it has proved to be a highly nontrivial problem
\cite{sub}. 

To define what is meant by the number of ends of a finitely generated
group take \(S \subseteq G\) a finite generating set of \(G\) and let
\( \Gamma (G, S)\) be the Cayley graph of \(G\) with respect to \(S\).
Then the number of ends is \(e(\Gamma (G, S))\) (\(e\) stands for edges)
which does not depend on the choice of a finite generating set \(S\)
of \(G\) hence it is well-defined.  Stallings' theorem about ends of
groups states that a finitely generated group \(G\) has more than one
end if and only if the group \(G\) admits a nontrivial decomposition
as an amalgamated free product or an HNN extension over a finite
subgroup \cite{end}. A word hyperbolic group roughly speaking, is a
finitely generated group equipped with a word metric satisfying
certain properties characteristic of hyperbolic geometry.

The famous Gromov's question has been much speculated about but it is
still very much open even for very concrete groups. It is not even
quite clear which answer to expect.  One of the few classes of groups
the answer is know to is in the case of Coxeter groups and some Artin
groups, where it is true \cite{coxeter}. It is not even know for
one-relator groups \(G_n(w) = \frac{F_n}{\langle \langle
  w\rangle\rangle}\) where \(w\) is an element of a free group of rank
\(n\), \(F_n\) that is not a proper power.

In this paper we focus on doubles \(D_n(w)= F_n *_{\langle w\rangle}
F_n\) where \(F_n\) is a free group of rank \(n\) and \(w \in F_n\).
The useful recent reduction of the question in the case of the doubles
is:

\begin{thm}[Gordon, Wilton]\cite{simple}
  \label{th:simple}
  Let \(w \in F_n\). If \(G_n(w)\) has an index-k subgroup \(G'\) with
  \(\beta_1 (G') > 1 + k(n - 2)\) then the double \(D_n(w)= F_n
  *_{\langle w\rangle} F_n\) contains a surface subgroup.
\end{thm}

The above result allowed Gordon and Wilton to exhibit several
infinite families of new examples of doubles with surface subgroups
\cite{simple}.  This result reduces the Gromov's question for doubles
to virtual homology. The only difficulty is that general approach to
computing the virtual homology is not fully developed. But for each
particular group \(\beta_1\) of a subgroup can be attempted to be
calculated using a computer. This is the approach taken to gather
evidence for Gromov's question.

We will be using this to investigate mainly doubles with \(n=3\). It will
be also shown how this method works for \(n=4\). The
Question~\ref{co:qu} was already studied with success  by Button in the case of
\(n=2\) \cite{2rel} using similar methods.


\section{Algorithm}

We will be looking at groups of the form \(G_3(w)=\langle a,b,c \mid w(a,b ,c)\rangle\)
where \(w(a,b,c)\) is a cyclically reduced word in three
letters of length up to 18. The reason 18 is chosen is that the longer
the word is the more computational time is needed. Cyclically reduced means
that cyclic
permutations are reduced. Reduced simply means that all obvious
cancellations like \(a^{-1}\) followed by \(a\) are made.

In general there is no algorithm to decided weather a group is
hyperbolic of not. But in the one relator setting there is a number of
theorems we will need to use later on.  

The trivial
corollary to the above Theorem~\ref{th:simple} that we will be using
in the remainder of the paper is:
\begin{cor}
  \label{co:main}
  Let \(w \in F_3\). If \(G_3(w)\) has an index-k subgroup \(G'\) with
  \(\beta_1 (G') > n+1\) then the double \(D_3(w)= F_3*_{\langle
    w\rangle} F_3\) contains a surface subgroup.
\end{cor}

A randomly chosen finitely presented group is almost surely
word-hyperbolic with an appropriate definition of ``almost surely''.
That is why initially we did not included any checks weather the
group is hyperbolic or not. Double groups are one-ended if \(w\) is not in a
proper free factor of \(F_n\)  \cite{simple}.

The algorithm is as follows:
\begin{enumerate}
\item Generate a random word \(w(a,b,c)\). It is done by choosing randomly
  18 characters from \(a, a^{-1}, b, b^{-1}, c, c^{-1}\) and then
  cyclically reduce the word. Since we do not know how much cancellation will take place we
  only know that the resulting word will be of length smaller then 18
  typically around 14.
\item Calculate the  index \(i=1\) subgroups  for \(G_3(w)=\langle a,b,c
  \mid w(a,b ,c)\rangle\).
\item Checking for the condition in
  Corollary \ref{co:main} for each subgroup. So for each subgroup  we calculate the first Betti number.
  In other words the abelinisation of the subgroup is calculated and
  the first betti number is the rank of it.
\item If the condition is  satisfied to step 6.
\item If the condition is not satisfied go back to step 2 and increase
 \( i\) by \(1\). Do this until \(i<10\) then move to the next step. The reason why the index is chosen to be 10 is that it is
  the highest average computer will calculate in reasonable time for a
  generic group.
\item Record the result weather the condition is satisfied for all \(i\) and go to step 1.  The output is \(w(a,b,c)\) and either the
  program found that the condition is satisfied and if so at which
  index or that it failed. 
\item Then calculate the number of successes and fails over the
   number of groups tried.
\end{enumerate}

This algorithm was implemented in MAGMA software for symbolic
calculations. The main limitation of this method is the speed of the
computer.

\section{Analysis of results}

It became clear after running the program that \(F_2\) appears quite
often which produces double that are not hyperbolic (and not
one-ended). To filter it out we used three methods.

The first one is linked to the Nielsen's moves. Let \(G\) be a group and
let \(M = (g_1,... , g_n) \in
G_n\) be an \(n\)-tuple of elements of \(G\).
The following moves are called elementary
Nielsen moves on \(M\), for  generators \(g_i\), \(1< i <n\) :
\begin{enumerate}
\item For some \(i\), replace \(g_i\) by \(g^{-1}\) in \(M\).
\item For some \(i \ne j\), \(1 <i, j <n\) replace \(g_i\) by \(g_ig_j\) in \(M\).
\item For some \(i \ne j\), \(1 <i\), \(j <n\) interchange \(g_i\) and \(g_j\)
  in \(M\).
\end{enumerate}
We say that two n-tuples are Nielsen equivalent  if there is a chain of
elementary Nielsen moves
which transforms one into another. In fact if they are Nielsen
equivalent if and only if they generate the same group. So if \(w\) has only one of \(a, a^{-1}, b, b^{-1}, c, c^{-1}\) then
\(G_3(w)\) is Nielsen equivalent to \(F_2\). Hence this is the first
thing to check for, since it is the least computationally expensive.

Secondly there is a function in MAGMA which looks for isomorphisms
between groups. So the next job is to find  isomorphisms of
\(G_3(w)\) and \(F_2\) which we do with parameter 9. The parameter in
the function isomorphism indicated how hard it looks in the sense the
higher the parameter the longer it will try to look for before giving
up. 

 Finally we can gather evidence that the group is \(F_2\) or at
least disprove that it is not by looking at the number of subgroups of
a all index (up to conjugacy the way it is counted in MAGMA). If the
group has the same number of subgroups for all indexes up to \(9\) it
is likely to be either \(F_2\) or it is indistinguishable (by looking
at subgroups) from it.

The second class of groups which will not satisfy the
condition  in Corollary \ref{co:main}  is:
when it can be written as \(G=F_1*H\) where \(H=BS(n,m)\)
Baumslag--Solitar groups or very close to them like \(\langle a,b^p \mid
b^{-p}a^nb^p=a^{m\pm1}\rangle\). The later are  due to Higman \cite{hig} called
Baumslag--Brunner--Gersten in \cite{2rel}. Baumslag--Solitar groups
are of the form:
\begin{displaymath}
  BS(n,m) \cong \langle a,b \mid b^{-1}a^nb=a^m\rangle.
\end{displaymath}
We have \(\beta_1(G')=i+1\) for all subgroups \(G'\) of \(G\) with
index \(i\) see section 5.5 in \cite{pro} for the proof.

\begin{prop}
  \label{pr:BBG}
  The linked one relator groups which do not give a surface subgroups
  for the doubles in light of Corollary \ref{co:main} are either
  \begin{enumerate}
  \item The free group on two generators \(F_2\).
  \item The groups of the form \(G=F_1*H\) where \(H=BS(n,m)\)
    Baumslag--Solitar or Baumslag--Brunner--Gersten groups.
  \end{enumerate}
\end{prop}

Note that for all subgroups \(F'\) of a free group on two generators
\(F_2\) we have that \(\beta_1(F')=i+1\). 

It appears that those can be very not trivial to spot even using a
computer. What is easier is to prove that the group is not of the a
certain form. In the next section we display four examples where our
computer did not find any surface subgroups but which are not of the
above kind.

The above algorithm was run for 1000 random groups showing that it
either contains a surface subgroup or are of the above form for \(96
\%\) of all group.

\section{Open Questions}

The relator of  four \(G_3(w)\) groups which are not of the above form are:

\begin{itemize}
\item \(ba^{-1}c^{-1}b^{-1}ab^{-2}a^{-1}b^{-1}a^{-1}c\)
\item \(b^{-1}ab^{-3}c^{-2}b^{-2}ca^{-1}\)
\item \(a c^{-1}  a  c^{-1}  a  c^2  b^2  c  a\)
\item \(a^{-1}  b  c^2  a^{-1}  c  a^{-1}  b  a^3  b^{-1}\)
\end{itemize}

One way to see that those relators do not give rise to \(F_2\) is to
compare the number of subgroups of index say \(9\) which is different
for each one. This actually indicates that no two of the above groups
are isomorphic to each other. To see that it is not \(G=F_1*H\) where
\(H=BS(n,m)\) we prove that all four of the above groups are word
hyperbolic.  To do that we use the paper \cite{hyper}, which has a nice criteria
in the case of one relator groups. The approach works if we can find a
presentation which has one letter appearing no more then three times,
which we have in all of the above. Then it is simply the matter of an
easy check. Note that hyperbolic groups cannot contain a
Baumslag--Solitar groups as a subgroup, this implies that the above
groups are not of the form as in the Proposition~\ref{pr:BBG}.

The above groups are intriguing: could it be the case that they are
decomposable but not with Baumslag--Solitar groups? If this is the
case then  \(G = F_1* H\), where \(H\) will be two generator one
relator group and the relator could be either of height 1 or not. One
relator groups are well studied and all of the counter-examples seem
to come from height 1 relators. Without loss of generality height one
word is:
\begin{displaymath}
  w=ba^{i_1}ba^{i_2} ... ba^{i_{l-1}}ba^{i_l}
\end{displaymath}
If it is not height 1 we cannot say anything about that at the moment.
But if \(H\) is then there is a theorem of J.~Button in \cite{hight}
which says that those groups are either large or are indistinguishable
from Baumslag-Solitar groups from looking at subgroups. If \(H\) is
large then certainly the \( \beta_1 > 1 + i\) for some \(i\). And also
the above groups do not have the same number of subgroups as any
Baumslag--Solitar groups. 

The way to see that is to note that we can work out the \(n - m\) from
the abelianisation and bound \(n + m < 16\) by the fact
that Nielsen's moves preserve the highest powers. Then there is only a
few possible Baumslag--Solitar groups to check, and none of them work
for any of the above groups. In fact the number of subgroups is
strictly in between that of \(F_2\) and Baumslag-Solitar groups. Hence
if it is decomposable then \(H\) is not of height \(1\).

It maybe it  might be the case that a higher
index is needed to detect the required
property. Or do there exist doubles which have a surface subgroup but
this is not detectable by Theorem~\ref{th:simple} for arbitrary index?  The
property that the above groups seem to share is very little torsion in
the abelianisation of subgroups. Also up to index \(9\) there is no
abelianisation of a subgroup which has the repeating torsion, which
could have been used to try the method described in Section 7.

\section{Decomposable into the two generators one relator group and a free group}

In this section we will be testing the groups of the form \(F(b, c
\mid w(b,c))* F_1\) where \(F_1\) is the free one generator group.
It is interesting to see  for which  \(F(b, c
\mid w(b,c))\)  we cannot find a surface subgroup in the associated doubles.

Question~\ref{co:qu} does not apply in this senario (the above is not
one-ended) but it is still of interest to see how many of them
actually satisfy the above property. Using this approach we were able
to come up with examples of groups which will not terminate using the
below program but which nevertheless satisfy the condition in Result
\ref{co:main}. Since the group is a free product it is enough to study
\(F(a,b \mid w(b,c))\) which is much smaller and so the computer can
go to a much higher index. 

For example, with \(w(b,c)=c^{-1} * b^{-2} * c^2
* b^{-3} * c^{-2}\) the property is only detected at index \(13\). One
would need a very powerful computer to go that high for \(n=3\).
Furthermore for \(w(b,c)=c * b^2 * c * b * c * b^{-1} * c^{-2} * b\)
the index the property is detected is \(30\). It is not possible to
calculate all subgroups up to index 30 even with the most powerful
computer. So we use a trick that was used in \cite{2rel} by spotting
that at index \(15\) there is a subgroup with an abelianisation which
had three cyclic groups of the same order. So take this subgroup as
the group and repeat the process with it, where it works already at
index \(2\).

\section{Four Generators}

We also tried this method in the case of \(n=4\). The program below
dealt with about \(87 \%\) of the random double groups. The reason why less of
them is dealt with is that with more generators the algorithms become
more expensive and the index up to which it is possible to go is only
6. Also the index we might need to go up to
might be bigger.

The way the algorithm worked is as follows:
\begin{enumerate}
\item We pick a random relator in the same sense as in the three generator
case. 
\item Cyclically reduce it. 
\item Check if there is a letter which occurs only once or not at all.
\item If it occurs only once then it is isomorphic to \(F_3\) by the
  Nelson's moves same as in the \(n=3\) case.
\item If there is a letter absent then \(G_4=F_1*G_3\) so it can be
  recovered from the \(n=3\) case. It is important to see if it is
  decomposable since \(\beta_1(K*L)=\beta_1(K)*\beta_1(L)\).
\item If neither of the two happens we search for isomorphisms with
  \(F_3\) this time with parameter 7 (smaller one had to be chosen due
  to more time consuming search).  Then we follow exactly the same
  procedure as in the case of \(n=3\).
\end{enumerate}

\section{Acknowledgement}

I am very grateful to Dr Jack Button for suggesting this project and for
the very helpful discussion along the way.
\bibliography{abbrevmr,akisil,newgeometry,analyse,algebra,arare,aclifford,aphysics}

\appendix
\section{The code for three generator groups}
\LGinlinefalse\LGbegin\lgrinde
\L{\LB{\V{F}\<\V{a},_\V{b},_\V{c}\>:=_\V{FreeGroup}(\N{3});}}
\L{\LB{\V{F1}\<\V{a1},\V{b1}\>:=_\V{FreeGroup}(\N{2});}}
\L{\LB{\V{kon}:=\N{0};}}
\L{\LB{\V{free2}:=\N{0};}}
\L{\LB{\V{sub}:=[\N{1},_\N{3},_\N{7},_\N{26},_\N{97},_\N{624},_\N{4163},_\N{34470},_\N{314493}];}\Tab{51}{\C{}//number_of\CE{}}}
\L{\LB{\C{}//_subgroups_of}\Tab{17}{\2(F\_2\2)_for_indexes_up_to_9_to_check_against\CE{}}}
\L{\LB{}}
\L{\LB{\K{for}_\V{i\_1}_:=_\N{1}_\K{to}_\N{50}_\K{do}_\C{}//numbers_of_groups_checked\CE{}}}
\L{\LB{}\Tab{2}{\V{rel}:=\V{Id}(\V{F});}}
\L{\LB{}\Tab{2}{\V{c1}:=\N{0};}}
\L{\LB{}\Tab{2}{\V{c2}:=\N{0};}}
\L{\LB{}\Tab{2}{\V{c3}:=\N{0};}}
\L{\LB{}\Tab{2}{\K{for}_\V{i}_:=_\N{1}_\K{to}_\N{18}_\K{do}}}
\L{\LB{}\Tab{4}{\V{j}:=\V{Random}(\N{3},_\N{6});}}
\L{\LB{}\Tab{4}{\K{if}_\V{j}}\Tab{10}{\V{eq}_\N{1}_\K{then}_\V{rel}:=\V{rel}*\V{a};}}
\L{\LB{}\Tab{4}{\K{elif}_\V{j}_\V{eq}_\N{2}_\K{then}_\V{rel}:=\V{rel}*\V{a}\5\-\N{1};}}
\L{\LB{}\Tab{4}{\K{elif}_\V{j}_\V{eq}_\N{3}_\K{then}_\V{rel}:=\V{rel}*\V{b};}}
\L{\LB{}\Tab{4}{\K{elif}_\V{j}_\V{eq}_\N{4}_\K{then}_\V{rel}:=\V{rel}*\V{b}\5\-\N{1};}}
\L{\LB{}\Tab{4}{\K{elif}_\V{j}_\V{eq}_\N{5}_\K{then}_\V{rel}:=\V{rel}*\V{c};}}
\L{\LB{}\Tab{4}{\K{else}_\V{rel}:=\V{rel}*\V{c}\5\-\N{1};}}
\L{\LB{}\Tab{4}{\K{end}_\K{if};}}
\L{\LB{}\Tab{2}{\K{end}_\K{for};}}
\L{\LB{}\Tab{2}{}}
\endlgrinde\LGend

\LGinlinefalse\LGbegin\lgrinde
\L{\LB{}\Tab{2}{\K{for}_\V{i}:=\N{0}_\K{to}_\V{\#rel}_\K{do}}\Tab{23}{\C{}//cyclically_reducing\CE{}}}
\L{\LB{}\Tab{4}{\V{l1}:=\V{LeadingGenerator}(\V{rel});}}
\L{\LB{}\Tab{4}{\V{rel1}:=\V{rel}*\V{l1};}}
\L{\LB{}\Tab{4}{\K{if}_\V{\#rel}_\V{gt}_\V{\#rel1}_\K{then}_\V{rel}:=\V{l1}\5\-\N{1}*\V{rel}*\V{l1};}}
\L{\LB{}\Tab{4}{\K{else}_\K{break};}}
\L{\LB{}\Tab{4}{\K{end}_\K{if};}}
\L{\LB{}\Tab{2}{\K{end}_\K{for};}}
\L{\LB{}\Tab{2}{}}
\L{\LB{}\Tab{2}{\V{seq}:=\V{Eltseq}(\V{rel});}}
\L{\LB{}\Tab{2}{}}
\L{\LB{}\Tab{2}{\V{k}:=\N{1};}}
\L{\LB{}\Tab{2}{\K{for}_\V{i}_:=_\N{1}_\K{to}_\V{\#rel}_\K{do}}\Tab{25}{\C{}//counting_the_number_of_each_relator\CE{}}}
\L{\LB{}\Tab{4}{\K{if}_\V{seq}[\V{i}]}\Tab{15}{\V{eq}_\N{1}_\K{then}_\V{c1}:=\V{c1}+\N{1};}}
\L{\LB{}\Tab{4}{\K{elif}}\Tab{10}{\V{seq}[\V{i}]_\V{eq}_\-\N{1}_\K{then}_\V{c1}:=\V{c1}+\N{1};}}
\L{\LB{}\Tab{4}{\K{elif}_\V{seq}[\V{i}]_\V{eq}_\N{2}_\K{then}_\V{c2}:=\V{c2}+\N{1};}}
\L{\LB{}\Tab{4}{\K{elif}}\Tab{10}{\V{seq}[\V{i}]\V{eq}_\-\N{2}_\K{then}_\V{c2}:=\V{c2}+\N{1};}}
\L{\LB{}\Tab{4}{\K{elif}_\V{seq}[\V{i}]_\V{eq}_\N{3}_\K{then}_\V{c3}:=\V{c3}+\N{1};}}
\L{\LB{}\Tab{4}{\K{else}_\V{c3}:=\V{c3}+\N{1};}}
\L{\LB{}\Tab{4}{\K{end}_\K{if};}}
\L{\LB{}\Tab{2}{\K{end}_\K{for};}}
\L{\LB{}\Tab{2}{\K{if}_\V{c1}}\Tab{9}{\V{eq}_\N{1}_\K{then}_\V{k}:=\N{2};_\K{print}_\S{}\3Isomorphic_to_F2_trivially\3\SE{};_\V{free2}:=\V{free2}+\N{1};}}
\L{\LB{}\Tab{2}{\K{elif}_\V{c3}_\V{eq}_\N{1}_\K{then}_\V{k}:=\N{2};_\K{print}_\S{}\3Isomorphic_to_F2_trivially\3\SE{};_\V{free2}:=\V{free2}+\N{1};}}
\L{\LB{}\Tab{2}{\K{elif}_\V{c2}_\V{eq}_\N{1}_\K{then}_\V{k}:=\N{2};_\K{print}_\S{}\3Isomorphic_to_F2_trivially\3\SE{};_\V{free2}:=\V{free2}+\N{1};}}
\L{\LB{}\Tab{2}{\K{end}_\K{if};}}
\L{\LB{}\Tab{2}{\V{rel};}}
\L{\LB{}\Tab{2}{\V{G}_\<\V{e},_\V{f},_\V{g}\>_:=_\V{quo}\<\V{F}_\|_\V{rel}_\>;}}
\L{\LB{}\Tab{2}{\V{ab}:=\N{0};}}
\L{\LB{}\Tab{2}{\V{nu}:=\N{0};}}
\L{\LB{}\Tab{2}{\V{sB}:=\N{0};}}
\L{\LB{}\Tab{2}{\K{if}_\V{k}_\V{eq}_\N{1}_\K{then}}}
\L{\LB{}\Tab{4}{\V{isiso},_\V{f1},_\V{f2}_:=_\V{SearchForIsomorphism}(\V{G},\V{F1},\N{9});}}
\L{\LB{}\Tab{4}{\V{isiso};}}
\L{\LB{}\Tab{4}{\K{if}_\V{isiso}_\K{then}_\V{k}:=\N{2};_\K{print}_\S{}\3Isomorphic_to_F2\3\SE{};_\V{free2}:=\V{free2}+\N{1};}}
\L{\LB{}\Tab{4}{\K{end}_\K{if};}}
\L{\LB{}\Tab{2}{\K{end}_\K{if};}}
\L{\LB{}\Tab{2}{}}
\endlgrinde\LGend
\vfill
\newpage
\LGinlinefalse\LGbegin\lgrinde
\L{\LB{}\Tab{2}{\K{for}_\V{i}_:=_\N{1}_\K{to}}\Tab{17}{\N{9}}\Tab{20}{\K{do}_\C{}//the_index_up_to_which_it_is_going_up\CE{}}}
\L{\LB{}\Tab{4}{\K{if}_\V{k}_\V{eq}_\N{2}_\K{then}_\K{break};}}
\L{\LB{}\Tab{4}{\K{end}_\K{if};}}
\L{\LB{}\Tab{4}{\V{t}:=\V{LowIndexSubgroups}(\V{G},_\<\V{i},_\V{i}\>);}}
\L{\LB{}\Tab{4}{\K{if}_\V{\#t}_\V{ne}_\V{sub}[\V{i}]_\K{then}_\V{sB}:=\N{1};}}
\L{\LB{}\Tab{4}{\K{end}_\K{if};}}
\L{\LB{}\Tab{4}{\K{for}_\V{j}:=_\N{1}_\K{to}_\V{\#t}_\K{do}}}
\L{\LB{}\Tab{6}{\V{l}:=\V{AQInvariants}(\V{t}[\V{j}]);}}
\L{\LB{}\Tab{2}{}}
\L{\LB{}\Tab{6}{\V{con}:=\N{0};_\C{}//calculating_the_number_of_zero{'}s_in_the_abelinisation_\CE{}}}
\L{\LB{}\Tab{6}{\K{for}_\V{m}:=\N{1}_\K{to}_\V{\#l}_\K{do}}}
\L{\LB{}\Tab{8}{\K{if}_\N{1}_\V{gt}_\V{l}[\V{m}]_\K{then}_\V{con}:=\V{con}+\N{1};}}
\L{\LB{}\Tab{8}{\K{else}_\V{ab}:=\N{1};}}
\L{\LB{}\Tab{8}{\K{end}_\K{if};}}
\L{\LB{}\Tab{6}{\K{end}_\K{for};}}
\L{\LB{}\Tab{2}{}}
\L{\LB{}\Tab{6}{\K{if}_\V{con}_\V{gt}_\V{i}+\N{1}_\K{then}_\K{print}_\V{i};}\Tab{35}{\V{k}:=\N{2};_\V{kon}:=\V{kon}+\N{1};_\C{}//checking_condition\CE{}}}
\L{\LB{}\Tab{6}{\K{end}_\K{if};}}
\L{\LB{}\Tab{6}{\K{if}_\V{con}_\V{ne}_\V{i}+\N{1}_\K{then}_\V{nu}:=\N{1};}}
\L{\LB{}\Tab{6}{\K{end}_\K{if};}}
\L{\LB{}\Tab{2}{}}
\L{\LB{}\Tab{6}{\K{if}_\V{k}_\V{eq}_\N{2}_\K{then}_\K{break};}}
\L{\LB{}\Tab{6}{\K{end}_\K{if};}}
\L{\LB{}\Tab{4}{\K{end}_\K{for};}}
\L{\LB{}\Tab{2}{}}
\L{\LB{}\Tab{4}{\K{if}_\V{k}_\V{eq}_\N{2}_\K{then}_\K{break};}}
\L{\LB{}\Tab{4}{\K{end}_\K{if};}}
\L{\LB{}\Tab{2}{\K{end}_\K{for};}}
\L{\LB{}\Tab{2}{\K{if}_\V{k}_\V{eq}_\N{1}_\K{then}_\K{print}_\S{}\3Did_not_find_surface_subgroups\3\SE{};}}
\L{\LB{}\Tab{2}{\K{end}_\K{if};}}
\L{\LB{}\Tab{2}{\K{if}_\V{ab}_\V{eq}_\N{0}}\Tab{14}{\K{and}_\V{nu}_\V{eq}_\N{0}_\K{and}_\V{sB}_\V{eq}_\N{0}_\K{then}_\K{print}_\S{}\3Looks_like_F2\3\SE{};}}
\L{\LB{}\Tab{2}{\K{end}_\K{if};}}
\L{\LB{\K{end}_\K{for};}}
\L{\LB{\K{print}_\S{}\3Free_2\3\SE{};}}
\L{\LB{\V{free2};}}
\L{\LB{\K{print}_\S{}\3Done\3\SE{};}}
\L{\LB{\V{kon};}}
\endlgrinde\LGend
\vfill
\section{The code for four generator groups}

\LGinlinefalse\LGbegin\lgrinde
\L{\LB{\V{F}\<\V{a},_\V{b},_\V{c},_\V{d}\>:=_\V{FreeGroup}(\N{4});}}
\L{\LB{\V{F1}\<\V{a1},\V{b1},_\V{c1}\>:=_\V{FreeGroup}(\N{3});}}
\L{\LB{\V{kon}:=\N{0};}}
\L{\LB{\V{free3}:=\N{0};}}
\L{\LB{\V{free2}:=\N{0};}}
\L{\LB{\V{sub}:=[\N{1},_\N{7},_\N{41},_\N{604},_\N{13753},_\N{504243}];_\C{}//number_of\CE{}}}
\L{\LB{\C{}//_subgroups_of}\Tab{17}{\2(F\_3\2)_for_indexes_up_to_6}\Tab{46}{to_check_against\CE{}}}
\L{\LB{}}
\L{\LB{\K{for}_\V{i\_1}_:=_\N{1}_\K{to}_\N{50}}\Tab{20}{\K{do}}}
\L{\LB{}\Tab{2}{\V{rel}:=\V{Id}(\V{F});}}
\L{\LB{}\Tab{2}{\V{c1}:=\N{0};}}
\L{\LB{}\Tab{2}{\V{c2}:=\N{0};}}
\L{\LB{}\Tab{2}{\V{c3}:=\N{0};}}
\L{\LB{}\Tab{2}{\V{c4}:=\N{0};}}
\L{\LB{}\Tab{2}{\K{for}_\V{i}_:=_\N{1}_\K{to}_\N{14}_\K{do}}}
\L{\LB{}\Tab{4}{\V{j}:=\V{Random}(\-\N{1},_\N{6});}}
\L{\LB{}\Tab{4}{\K{if}_\V{j}}\Tab{10}{\V{eq}_\N{1}_\K{then}_\V{rel}:=\V{rel}*\V{a};}}
\L{\LB{}\Tab{4}{\K{elif}_\V{j}_\V{eq}_\-\N{1}_\K{then}_\V{rel}:=\V{rel}*\V{d}\5\-\N{1};}}
\L{\LB{}\Tab{4}{\K{elif}_\V{j}_\V{eq}_\N{0}_\K{then}_\V{rel}:=\V{rel}*\V{d};}}
\L{\LB{}\Tab{4}{\K{elif}_\V{j}_\V{eq}_\N{2}_\K{then}_\V{rel}:=\V{rel}*\V{a}\5\-\N{1};}}
\L{\LB{}\Tab{4}{\K{elif}_\V{j}_\V{eq}_\N{3}_\K{then}_\V{rel}:=\V{rel}*\V{b};}}
\L{\LB{}\Tab{4}{\K{elif}_\V{j}_\V{eq}_\N{4}_\K{then}_\V{rel}:=\V{rel}*\V{b}\5\-\N{1};}}
\L{\LB{}\Tab{4}{\K{elif}_\V{j}_\V{eq}_\N{5}_\K{then}_\V{rel}:=\V{rel}*\V{c};}}
\L{\LB{}\Tab{4}{\K{else}_\V{rel}:=\V{rel}*\V{c}\5\-\N{1};}}
\L{\LB{}\Tab{4}{\K{end}_\K{if};}}
\L{\LB{}\Tab{2}{\K{end}_\K{for};}}
\L{\LB{}\Tab{2}{}}
\L{\LB{}\Tab{2}{\K{for}_\V{i}:=\N{0}_\K{to}_\V{\#rel}_\K{do}}}
\L{\LB{}\Tab{4}{\V{l1}:=\V{LeadingGenerator}(\V{rel});}}
\L{\LB{}\Tab{4}{\V{rel1}:=\V{rel}*\V{l1};}}
\L{\LB{}\Tab{4}{\K{if}_\V{\#rel}_\V{gt}_\V{\#rel1}_\K{then}_\V{rel}:=\V{l1}\5\-\N{1}*\V{rel}*\V{l1};}}
\L{\LB{}\Tab{4}{\K{else}_\K{break};}}
\L{\LB{}\Tab{4}{\K{end}_\K{if};}}
\L{\LB{}\Tab{2}{\K{end}_\K{for};}}
\L{\LB{}\Tab{2}{}}
\L{\LB{}\Tab{2}{\V{seq}:=\V{Eltseq}(\V{rel});}}
\L{\LB{}\Tab{2}{}}
\L{\LB{}\Tab{2}{\V{k}:=\N{1};}}
\L{\LB{}\Tab{2}{\K{for}_\V{i}_:=_\N{1}_\K{to}_\V{\#rel}_\K{do}}}
\L{\LB{}\Tab{4}{\K{if}_\V{seq}[\V{i}]}\Tab{15}{\V{eq}_\N{1}_\K{then}_\V{c1}:=\V{c1}+\N{1};}}
\L{\LB{}\Tab{4}{\K{elif}}\Tab{10}{\V{seq}[\V{i}]_\V{eq}_\-\N{1}_\K{then}_\V{c1}:=\V{c1}+\N{1};}}
\L{\LB{}\Tab{4}{\K{elif}_\V{seq}[\V{i}]_\V{eq}_\N{2}_\K{then}_\V{c2}:=\V{c2}+\N{1};}}
\L{\LB{}\Tab{4}{\K{elif}}\Tab{10}{\V{seq}[\V{i}]\V{eq}_\-\N{2}_\K{then}_\V{c2}:=\V{c2}+\N{1};}}
\L{\LB{}\Tab{4}{\K{elif}_\V{seq}[\V{i}]_\V{eq}_\N{3}_\K{then}_\V{c3}:=\V{c3}+\N{1};}}
\L{\LB{}\Tab{4}{\K{elif}_\V{seq}[\V{i}]_\V{eq}_\N{4}_\K{then}_\V{c4}:=\V{c4}+\N{1};}}
\L{\LB{}\Tab{4}{\K{elif}_\V{seq}[\V{i}]_\V{eq}_\-\N{4}_\K{then}_\V{c4}:=\V{c4}+\N{1};}}
\L{\LB{}\Tab{4}{\K{else}_\V{c3}:=\V{c3}+\N{1};}}
\L{\LB{}\Tab{4}{\K{end}_\K{if};}}
\L{\LB{}\Tab{2}{\K{end}_\K{for};}}
\endlgrinde\LGend
\vfill
\newpage
\LGinlinefalse\LGbegin\lgrinde
\L{\LB{}\Tab{2}{\K{if}_\V{c1}}\Tab{9}{\V{eq}_\N{1}_\K{then}_\V{k}:=\N{2};_\K{print}_\S{}\3Isomorphic_to_F3_trivially\3\SE{};_\V{free3}:=\V{free3}+\N{1};}}
\L{\LB{}\Tab{2}{\K{elif}_\V{c3}_\V{eq}_\N{1}_\K{then}_\V{k}:=\N{2};_\K{print}_\S{}\3Isomorphic_to_F3_trivially\3\SE{};_\V{free3}:=\V{free3}+\N{1};}}
\L{\LB{}\Tab{2}{\K{elif}_\V{c2}_\V{eq}_\N{1}_\K{then}_\V{k}:=\N{2};_\K{print}_\S{}\3Isomorphic_to_F3_trivially\3\SE{};_\V{free3}:=\V{free3}+\N{1};}}
\L{\LB{}\Tab{2}{\K{elif}_\V{c4}_\V{eq}_\N{1}_\K{then}_\V{k}:=\N{2};_\K{print}_\S{}\3Isomorphic_to_F3_trivially\3\SE{};_\V{free3}:=\V{free3}+\N{1};}}
\L{\LB{}\Tab{2}{\K{elif}_\V{c1}}\Tab{11}{\V{eq}_\N{0}_\K{then}_\V{k}:=\N{2};_\K{print}_\S{}\3Back_to_3_generator_case\3\SE{};_\V{free2}:=\V{free2}+\N{1};}}
\L{\LB{}\Tab{2}{\K{elif}_\V{c3}_\V{eq}_\N{0}_\K{then}_\V{k}:=\N{2};_\K{print}_\S{}\3Back_to_3_generator_case\3\SE{};_\V{free2}:=\V{free2}+\N{1};}}
\L{\LB{}\Tab{2}{\K{elif}_\V{c2}_\V{eq}_\N{0}_\K{then}_\V{k}:=\N{2};_\K{print}_\S{}\3Back_to_3_generator_case\3\SE{};_\V{free2}:=\V{free2}+\N{1};}}
\L{\LB{}\Tab{2}{\K{elif}_\V{c4}_\V{eq}_\N{0}_\K{then}_\V{k}:=\N{2};_\K{print}_\S{}\3Back_to_3_generator_case\3\SE{};_\V{free2}:=\V{free2}+\N{1};}}
\L{\LB{}\Tab{2}{\K{end}_\K{if};}}
\L{\LB{}\Tab{2}{\V{rel};}}
\L{\LB{}\Tab{2}{\V{G}_\<\V{e},_\V{f},_\V{g}\>_:=_\V{quo}\<\V{F}_\|_\V{rel}_\>;}}
\L{\LB{}\Tab{2}{\V{ab}:=\N{0};}}
\L{\LB{}\Tab{2}{\V{nu}:=\N{0};}}
\L{\LB{}\Tab{2}{\V{sB}:=\N{0};}}
\L{\LB{}\Tab{2}{\K{if}_\V{k}_\V{eq}_\N{1}_\K{then}}}
\L{\LB{}\Tab{4}{\V{isiso},_\V{f1},_\V{f2}_:=_\V{SearchForIsomorphism}(\V{G},\V{F1},\N{7});}}
\L{\LB{}\Tab{4}{\V{isiso};}}
\L{\LB{}\Tab{4}{\K{if}_\V{isiso}_\K{then}_\V{k}:=\N{2};_\K{print}_\S{}\3Isomorphic_to_F2\3\SE{};_\V{free3}:=\V{free3}+\N{1};}}
\L{\LB{}\Tab{4}{\K{end}_\K{if};}}
\L{\LB{}\Tab{2}{\K{end}_\K{if};}}
\L{\LB{}\Tab{2}{}}
\L{\LB{}\Tab{2}{\K{for}_\V{i}_:=_\N{1}_\K{to}}\Tab{17}{\N{6}}\Tab{20}{\K{do}}}
\L{\LB{}\Tab{4}{\K{if}_\V{k}_\V{eq}_\N{2}_\K{then}_\K{break};}}
\L{\LB{}\Tab{4}{\K{end}_\K{if};}}
\L{\LB{}\Tab{4}{\V{t}:=\V{LowIndexSubgroups}(\V{G},_\<\V{i},_\V{i}\>);}}
\L{\LB{}\Tab{4}{\K{if}_\V{\#t}_\V{ne}_\V{sub}[\V{i}]_\K{then}_\V{sB}:=\N{1};}}
\L{\LB{}\Tab{4}{\K{end}_\K{if};}}
\L{\LB{}\Tab{4}{\K{for}_\V{j}:=_\N{1}_\K{to}_\V{\#t}_\K{do}}}
\L{\LB{}\Tab{4}{\V{l}:=\V{AQInvariants}(\V{t}[\V{j}]);}}
\L{\LB{}\Tab{2}{}}
\L{\LB{}\Tab{4}{\V{con}:=\N{0};}}
\L{\LB{}\Tab{4}{\K{for}_\V{m}:=\N{1}_\K{to}_\V{\#l}_\K{do}}}
\L{\LB{}\Tab{6}{\K{if}_\N{1}_\V{gt}_\V{l}[\V{m}]_\K{then}_\V{con}:=\V{con}+\N{1};}}
\L{\LB{}\Tab{6}{\K{else}_\V{ab}:=\N{1};}}
\L{\LB{}\Tab{6}{\K{end}_\K{if};}}
\L{\LB{}\Tab{4}{\K{end}_\K{for};}}
\L{\LB{}\Tab{2}{}}
\L{\LB{}\Tab{4}{\K{if}_\V{con}_\V{gt}_\N{2}*\V{i}+\N{1}_\K{then}_\K{print}_\V{i};}\Tab{35}{\V{k}:=\N{2};_\V{kon}:=\V{kon}+\N{1};}}
\L{\LB{}\Tab{4}{\K{end}_\K{if};}}
\L{\LB{}\Tab{4}{\K{if}_\V{con}_\V{ne}_\V{i}+\N{1}_\K{then}_\V{nu}:=\N{1};}}
\L{\LB{}\Tab{4}{\K{end}_\K{if};}}
\L{\LB{}\Tab{2}{}}
\L{\LB{}\Tab{4}{\K{if}_\V{k}_\V{eq}_\N{2}_\K{then}_\K{break};}}
\L{\LB{}\Tab{4}{\K{end}_\K{if};}}
\L{\LB{}\Tab{4}{\K{end}_\K{for};}}
\L{\LB{}\Tab{4}{\K{if}_\V{k}_\V{eq}_\N{2}_\K{then}_\K{break};}}
\L{\LB{}\Tab{4}{\K{end}_\K{if};}}
\L{\LB{}\Tab{4}{\K{end}_\K{for};}}
\L{\LB{}\Tab{4}{\K{if}_\V{k}_\V{eq}_\N{1}_\K{then}_\K{print}_{`}{`}\V{Did}_\K{not}_\V{find}_\V{surface}_\V{subgroups}\S{}{'}{'}\SE{};}}
\L{\LB{}\Tab{4}{\K{end}_\K{if};}}
\L{\LB{}\Tab{4}{\K{if}_\V{ab}_\V{eq}_\N{0}_\K{and}_\V{nu}_\V{eq}_\N{0}_\K{and}_\V{sB}_\V{eq}_\N{0}_\K{then}_\K{print}_\S{}\3Looks_like_F3\3\SE{};}}
\L{\LB{}\Tab{4}{\K{end}_\K{if};}}
\L{\LB{\K{end}_\K{for};}}
\endlgrinde\LGend
\LGinlinefalse\LGbegin\lgrinde
\L{\LB{\K{print}_\S{}\3Free_3:\3\SE{};}}
\L{\LB{\V{free3};}}
\L{\LB{\K{print}_\S{}\3Back_to_3_generator_case:\3\SE{};}}
\L{\LB{\V{free2};}}
\L{\LB{\K{print}_\S{}\3Done:\3\SE{};}}
\L{\LB{\V{kon};}}
\endlgrinde\LGend

\end{document}